\documentclass[12pt,a4paper]{article}
\usepackage{amssymb}

\usepackage{amsmath}

\usepackage{amsfonts}

\numberwithin{equation}{section}
\newtheorem{thm}{Theorem}[section]
\newtheorem{prop}[thm]{Proposition}
\newtheorem{lem}[thm]{Lemma}
\newtheorem{cor}[thm]{Corollary}

\newtheorem{defn}[thm]{Definition}

\newfont{\gothic}{eufb10}

\title{Components of maximal dimension
in the Noether-Lefschetz locus
for Beilinson-Hodge cycles on open surfaces}
\author{M. Asakura and S. Saito}
\date{\empty}
\begin{document}
\maketitle

\vskip 20pt

\renewcommand{\labelenumi}{(\theenumi)}

\def\Th#1{\bigskip\noindent\textbf{Theorem #1}\; }
\def\Cor#1{\bigskip\noindent\textbf{Corollary #1}\; }
\def\Def#1{\bigskip\noindent\textbf{Definition #1}\; }
\def\Deflem#1{\bigskip\noindent\textbf{Definition-Lemma #1}\; }
\def\Lem#1{\bigskip\noindent\textbf{Lemma #1}\; }
\def\Rem#1{\bigskip\noindent\textbf{Remark #1}\; }
\def\Claim#1{\bigskip\noindent\textbf{Claim #1}\; }
\def\Prop#1{\bigskip\noindent\textbf{Proposition #1}\; }
\def\proof{\bigskip\noindent\textbf{Proof}\; }

\def\cO{{\mathcal{O}}}
\def\cX{{\mathcal{X}}}
\def\cZ{{\mathcal{Z}}}
\def\cY{{\mathcal{Y}}}
\def\cV{{\mathcal{V}}}
\def\cU{{\mathcal{U}}}

\def\bA{\Bbb A}

\def\codim{\mathrm{codim}}
\def\Im{\mathrm{Im}}
\def\Ker{\mathrm{Ker}}
\def\Coker{\mathrm{Coker}}
\def\Ext{\mathrm{Ext}}
\def\Hom{\mathrm{Hom}}
\def\PGL{\mathrm{PGL}}
\def\Pic{\mathrm{Pic}}

\def\mod{\; \mathrm{mod} \;}
\def\PG#1{\frak{S}_{#1}}

\def\tM{\tilde{M}}

\def\scs{\; : \;}
\def\NLl{\hbox{Noether-Lefschetz locus}}
\def\MNL{M_{NL}}
\def\SU{\Sigma(U)}
\def\SUt{\Sigma(U_t)}
\def\dlog#1{d\log #1}

\def\uc{\underline{c}}
\def\us{\underline{s}}
\def\uct{\underline{c}_t}
\def\Tij{T_{ij}}
\def\Tabc{T_{123}}
\def\Ta{T_{12}}
\def\Tb{T_{23}}
\def\Tc{T_{31}}
\def\Tijc{T_{ij}(\uc)}
\def\Tijcd{T_{i'j'}(\uc')}
\def\Tpqc{T_{(p,q)}(\uc)}
\def\Tpqcd{T_{(p',q')}(\uc')}
\def\Tspqc{T^{\sigma}_{(p,q)}(\uc)}
\def\Tspqcd{T^{\sigma}_{(p',q')}(\uc')}

\def\Dl{\Delta_{\lambda}}
\def\Dli{\Delta_{\lambda_i}}
\def\Dlw{\Delta_{\lambda,w}}
\def\Dlwq{\Delta_{\lambda,w,\Gamma}}
\def\Dlo{\Delta_{\lambda}^o}
\def\Dla{\Delta_{\lambda}^1}
\def\Dlb{\Delta_{\lambda}^1}
\def\Tl{T_{\lambda}}
\def\Wl{W_{\lambda}}
\def\Vl{V_{\lambda}}
\def\Rl{R_{\lambda}}
\def\Il#1{I_{\lambda}^{#1}}
\def\Ilt#1{I_{\lambda_t}^{#1}}
\def\Ill{I_{\lambda}}
\def\Pw#1{P_w^{#1}}
\def\Mwq{M_{w,\Gamma}}
\def\Mwe{M_{w}^{(e)}}
\def\Mwqe{M_{w,\Gamma}^{(e)}}
\def\Ge{G^{(e)}}
\def\Geq{G_{\Gamma}^{(e)}}
\def\Go{G^{(0)}}
\def\pie{\pi^{(e)}}
\def\Lq{L}

\def\RF#1{R_F^{#1}}
\def\RFt#1{R_{\Ft}^{#1}}
\def\JF#1{J_F^{#1}}
\def\JFt#1{J_{\Ft}^{#1}}
\def\IRF#1#2{<#1>_F^{#2}}

\def\cijX{c_{ij}(X)}
\def\cijXt{c_{ij}(X_t)}
\def\caX{c_{12}(X)}
\def\caXt{c_{12}(X_t)}
\def\cbX{c_{23}(X)}
\def\cbXt{c_{23}(X_t)}
\def\ccX{c_{31}(X)}
\def\ccXt{c_{31}(X_t)}

\def\Xt{X_t}
\def\Ut{U_t}
\def\Vt{U'_t}
\def\Zt{Z_t}
\def\Wt{W_t}
\def\Ft{F_t}
\def\Qt{\Gamma_t}
\def\Wm{W_{\mu}}
\def\at{\alpha_{\Ut}}
\def\bt{\beta_{\Ut}}
\def\tbt{\widetilde{\beta}_{\Ut}}
\def\wUt{\omega_{\Ut}}
\def\xt{\xi_{\Ut}}
\def\lt{\lambda(t)}
\def\lit{\lambda_{i}(t)}
\def\ait{a_i(t)}
\def\aat{a_1(t)}
\def\abt{a_2(t)}
\def\act{a_3(t)}
\def\aio{a_i(0)}

\def\aU{\alpha_{U}}
\def\wU{\omega_{U}}
\def\wF{\omega_{F}}
\def\eF{\eta_{F}}
\def\kF{\kappa_{F}}
\def\xU{\xi_{U}}
\def\xF{\xi_{F}}

\def\dt{\partial_{tame}}
\def\dd{\partial_{div}}
\def\du{\partial_U}
\def\dz{\partial_Z}

\def\regU#1{reg_U^{#1}}
\def\regUd#1{reg_{U'}^{#1}}
\def\regDU#1{reg_{D,U}^{#1}}
\def\regDUd#1{reg_{D,U'}^{#1}}
\def\regDWm#1{reg_{D,\Wm}^{#1}}
\def\regDW#1{reg_{D,W}^{#1}}
\def\regUt#1{reg_{\Ut}^{#1}}
\def\regD#1{reg_{D,#1}}
\def\oX{\rho_X}
\def\oXt{\rho_{\Xt}}

\def\rmapo#1{\overset{#1}{\longrightarrow}}
\def\gij{g_{ij}}
\def\nabb{\overline{\nabla}}

\def\WM#1{\Omega_M^{#1}}
\def\WMo#1{\Omega_{M,0}^{#1}}
\def\WXZ#1{\Omega_X^{#1}(\log Z)}
\def\TXZ{\Theta_X(\log Z)}
\def\WXZt#1{\Omega_{\Xt}^{#1}(\log \Zt)}
\def\TMo{T_0(M)}

\def\HO{H_{\cO}}
\def\HQ{H_{\Bbb Q}}
\def\HC{H_{\Bbb C}}

\def\HQcU#1#2{\HQ^{#1}(\cU/M)(#2)}
\def\HCcU#1{\HC^{#1}(\cU/M)}
\def\HOcU#1{\HO^{#1}(\cU/M)}
\def\HcU#1#2{H^{#1,#2}(\cU/M)}
\def\HU#1#2{H^{#1,#2}(U)}
\def\HUt#1#2{H^{#1,#2}(\Ut)}

\def\isom{\overset{\sim}{\longrightarrow}}
\def\dP{\overset{\vee}{\Bbb P}}
\def\alphab{\underline{\alpha}}
\def\pd#1#2{\frac{\partial #1}{\partial #2}}

\def\fat{f_1(t)}
\def\fao{f_1(0)}
\def\fbt{f_2(t)}
\def\wt{\omega(t)}
\def\kt{\kappa(t)}
\def\et{\eta(t)}
\def\snt{s_{\nu}(t)}
\def\sno{s_{\nu}(0)}

\section*{Introduction}

\medbreak

In the moduli space $M$ of smooth hypersurfaces of degree $d$ in $\Bbb P^3$
over $\Bbb C$, the locus of those surfaces that possess curves which are not
complete intersections of the given surface with another surface is called the
Noether-Lefschetz locus and denoted by $\MNL$.
One can show that $\MNL$ is the union of a countable number of closed
algebraic subsets of $M$. The classical theorem of Noether-Lefschetz affirms
that every component of $\MNL$ has positive codimension in $M$ when $d\geq 4$.
Note that the theorem is false if $d=3$ since a smooth cubic surface has
the Picard number $7$. Since the infinitesimal method in Hodge theory was
introduced in [CGGH] as a powerful tool to study $\MNL$, fascinating results
have been obtained concerning irreducible components of $\MNL$. First we have
the following.

\begin{thm}[\cite{G1}]\label{0-1}
For every irreducible component $T$ of $\MNL$,
$\codim(T) \geq d-3$.
\end{thm}

\medbreak

The basic idea of the proof of the result is to translate the problem in
the language of the infinitesimal variation of Hodge structures on a family of
hypersurfaces. Then, by the Poincar\'e residue representation of the cohomology
of a hypersurface, the result follows from the duality theorem for the Jacobian
ring associated to a hypersurface.
We note that the inequality is the best possible since the family
of hypersurfaces of degree $d\geq 3$ containing a line has codimension exactly
$d-3$. M.Green [G2] and C.Voisin [V] has shown the following striking theorem.

\begin{thm} \label{0-2}
If $d\geq 5$, the only irreducible component of $\MNL$ having
codimension $d-3$ is the family of surfaces of degree $d$ containing a line.
\end{thm}

\medbreak

In this paper we study an analog of the above problem in the context of
Beilinson's Hodge conjecture. For a quasi-projective smooth variety $U$ over
$\Bbb C$, the space of Beilinson-Hodge cycles is defined to be
\begin{displaymath}
F^0H^q(U,\Bbb Q(q)):=H^q(U,\Bbb Q(q))\cap F^qH^q(U,\Bbb C)
\end{displaymath}
where $H^q(U,\Bbb Q(q))$ is the singular cohomology with coefficient
$\Bbb Q(q)=(2\pi\sqrt{-1})^q\Bbb Q$ and $F^{\bullet}$ denotes the Hodge
filtration of the mixed Hodge structure on the singular cohomology defined by
Deligne [D].
Beilinson's conjecture claims the surjectivity of the regulator map
(cf. [Bl] and [Sch])
\begin{displaymath}
\regU q\scs CH^q(U,q)\otimes\Bbb Q \to F^0H^q(U,\Bbb Q(q))
\end{displaymath}
where $CH^q(U,q)$ is Bloch's higher Chow group.
Taking a smooth compactification $U\subset X$ with $Z=X\setminus U$,
a simple normal crossing divisor on $X$, we have the following formula for
the value of $\regU q$ on decomposable elements in $CH^q(U,q)$;
\begin{displaymath}
\regU q(\{g_1,\dots,g_q\})=
\dlog {g_1}\wedge\cdots \wedge \dlog {g_q}\in H^0(X,\WXZ q)=F^0H^q(U,\Bbb C),
\end{displaymath}
where $\{g_1,\dots,g_q\}\in CH^q(U,q)$ is the products of
$g_j\in CH^1(U,1)=\Gamma(U,\cO^*_{Zar})$.
Beilinson's conjecture is an analog of the Hodge conjecture which claims the
surjectivity of cycle class maps from Chow group to space of Hodge cycles on
projective smooth varieties. The conjecture is known to hold in case $q=1$
(cf. [J], Th.5.1.3) but open in general in case $q\geq 2$.
\medbreak

The main subject of study in this paper is the {\it Noether-Lefschetz locus for
Beilinson-Hodge cycles} on the complement of the union of a normal crossing
divisor in a surface in $\Bbb P^3$. Let $X,Y_1,\dots,Y_s \subset\Bbb P^3$ be
smooth surfaces intersecting transversally and put
\begin{equation}\label{(0-1)}
Z=\underset{1\leq j\leq s}{\cup} Z_j \text{ with } Z_j=X\cap Y_j,\;
U=X\setminus Z.
\end{equation}
Let $H^2(U,\Bbb Q(2))_{triv}$ be the image of the natural restriction map
\begin{displaymath}
H^2(\Bbb P^3 \setminus \underset{1\leq j\leq s}{\cup} Y_j,\Bbb Q(2)) \to
H^2(U,\Bbb Q(2)).
\end{displaymath}
One can show ([AS2], Lem.(2-1))
\begin{displaymath}
H^2(U,\Bbb Q(2))_{triv}=\regU 2 (CH^2(U,2)_{dec}),
\end{displaymath}
where $CH^2(U,2)_{dec}\subset CH^2(U,2)\otimes\Bbb Q$ is the so-called
decomposable part, the subspace generated by the image of the product map
$CH^1(U,1)\otimes CH^1(U,1)\to CH^2(U,2)$.
It implies that
\begin{displaymath}
H^2(U,\Bbb Q(2))_{triv}\subset\Im(\regU 2)\subset F^0H^2(U,\Bbb Q(2)).
\end{displaymath}
We define
$F^0H^2(U,\Bbb Q(2))_{prim}:=F^0H^2(U,\Bbb Q(2))/H^2(U,\Bbb Q(2))_{triv}$
called the space of primitive Beilinson-Hodge cycles.
\par

Now fix integers
$d\geq 1$ and $e_j\geq 1$ with $1\leq j\leq s$.
Let $M$ be the moduli space of
sets of hypersurfaces $(X,Y_1,\dots,Y_s)$ of degree $(d,e_1,\dots,e_s)$
which intersect transversally.
Let $(\cX,\cY_1,\dots,\cY_s)$ be the universal family over $M$ and put
$$\cZ=\cX \cap (\cup_{1\leq j\leq s} \cY_j), \quad \cU=\cX\setminus\cZ.$$
For $t\in M$ let $U_t\subset X_t\supset Z_t$ be the fibers of
$\cU \subset \cX \supset \cZ$.

\begin{defn}\label{d0-1}
The Noether-Lefschetz locus for Beilinson-Hodge cycles on
$\cU/M$ is
$$ \MNL=\{t\in M|\;  F^0H^2(U_t,\Bbb Q(2))_{prim}\not=0\}$$
\end{defn}

\medbreak

The analogy with the classical Noether-Lefschetz locus is explained as follows.
Instead of the map
\begin{displaymath}
H^2(\Bbb P^3 \setminus \underset{1\leq j\leq s}{\cup} Y_j,\Bbb Q(2))\to
H^2(U,\Bbb Q(2))\cap F^2H^2(U,\Bbb C)
\end{displaymath}
we consider
\begin{displaymath}
H^2(\Bbb P^3,\Bbb Q(1))\to H^2(X,\Bbb Q(1))\cap F^1H^2(X,\Bbb C).
\end{displaymath}
By noting that the space on the left hand side is generated by
the cohomology class of a hyperplane, and
that that on the right hand side is identified with $\Pic(X)\otimes\Bbb Q$,
the space defined in the same way as Def.\ref{d0-1} is
nothing but the classical Noether-Lefschetz locus.
\medbreak

One can show as before that $\MNL$ is the union of a countable number of
closed analytic subsets.
By the analogy a series of problems on $\MNL$ arise, the problems to show
the counterparts of Th.\ref{0-1} and Th.\ref{0-2} in the new context.
In [AS2] the following result is shown.

\begin{thm}\label{0-3}
Assume $d\geq 4$. For every irreducible component $T$ of $\MNL$,
\begin{displaymath}
\codim_M(T) \geq d+\min\{d,e_1,\dots,e_s\}-2.
\end{displaymath}
\end{thm}

\medbreak

We should note that the estimate in Th.\ref{0-3} is far from being optimal
to the contrary to the case of Th.\ref{0-1}. It is observed in the main theorem
of this paper (Th.(0-4) below) that the optimal estimate in some case is given
by a quadratic polynomial in $d$.
The basic strategy of the proof of Th.\ref{0-3} is the same as that of
Th.\ref{0-1}.
A new input is the theory of generalized Jacobian rings developed in [AS1],
which give an algebraic description of the cohomology of the open surface $U$.
In particular the duality theorem for such rings plays a crucial role.
\medbreak

In order to state the main result of this paper, which is considered a
counterpart of Th.\ref{0-2},
we need restrict ourselves to the special case that
$s=3$ and $e_j=1$ for $1\leq j\leq 3$. Let $P=\Bbb C[z_0,z_1,z_2,z_3]$ be the
homogeneuous coordinate ring of $\Bbb P^3$. For the rest of the paper we let
$M$ be the moduli space of hypersurface of degree $d$ in $\Bbb P^3$ which
transversally intersects $Y=\cup_{1\leq j\leq 3} Y_j$ with
$Y_j=\{z_j=0\} \subset \Bbb P^3$. Let $\cX/M$ be the universal family and
$\Xt$ be its fiber over $t\in M$ and put $\Ut=\Xt\setminus (\Xt \cap Y)$.
Let $\MNL\subset M$ be defined as in Def.\ref{d0-1}. In order to determine the
irreducible components of $\MNL$ of maximal dimension, we need introduce some
notations. For an integer $l>0$ let $P^l\subset P$
be the subspace of homogeneous polynomials of degree $l$.

\begin{defn}\label{d0-2}
For a pair $(p,q)$ of non-negative coprime integers such that $d=r(p+q)$ with
$r\in \Bbb Z$,
and $\uc=[c_\nu]_{1\leq\nu\leq r}=[c_1:\dots :c_r]\in \Bbb P^r(\Bbb C)$,
and $\sigma\in \PG 3$, the permutation group on $(1,2,3)$, we let
$\Tspqc\subset M$ be the subset of those surfaces defined by an equation
$$ F=w A +\prod_{1\leq \nu\leq r}
(c z_{\sigma(1)}^{p+q}-c_\nu z_{\sigma(2)}^p z_{\sigma(3)}^q)
\text{ for some } w\in P^1,\; A\in P^{d-1}, c\in \Bbb C^*.$$
\end{defn}
\medbreak

We will see the following facts (cf. \S1):
\begin{itemize}
\item[(1)]
$\Tspqc$ is smooth irreducible and $\codim_M(\Tspqc)=\binom{d+2}{2}-5$.
\item[(2)]
If $c_\nu$ is a root of unity for $1\leq \forall \nu\leq r$,
$\Tspqc\subset \MNL$.
\end{itemize}

We now state the main theorem in this paper.

\begin{thm}\label{0-4} Assume $d\geq 4$.
\begin{itemize}
\item[$(1)$]
For every irreducible component $T\subset \MNL$,
$\codim(T)\geq \binom{d+2}{2}-5$.
\item[$(2)$]
The equality holds if and only if
$T=\Tspqc$ for some $\sigma$, $(p,q)$ and $\uc=[c_\nu]_{1\leq\nu\leq r}$
such that $c_\nu$ is a root of unity for $1\leq \forall \nu\leq r$.
\item[$(3)$]
If $X$ is a general member of $\Tspqc$, $\regU 2$ is surjective so that
Beilinson's Hodge conjecture holds for $U=X\setminus (X\cap Y)$.
\end{itemize}
\end{thm}
\vskip 5pt

A key to the proof of Th.\ref{0-4} is a result by Otwinowska ([Ot], Th.2)
on the Hilbert function of graded algebras of dimension 0.
\medbreak

Finally we discuss an implication of Th.\ref{0-4} on the injectivity of the
regulator map. Let $X$ be a member of $M$.
We are interested in the regulator map to Deligne cohomology
\begin{displaymath}
\oX \scs CH^2(X,1)\otimes\Bbb Q \to H^3_D(X,\Bbb Q(2)),
\end{displaymath}
where $CH^2(X,1)$ is Bloch's higher Chow group defined to be the cohomology of
the complex
\begin{displaymath}
K_2(\Bbb C(X)) \overset{\dt}{\longrightarrow}
\bigoplus_{C\subset X} \Bbb C(C)^*
\overset{\dd}{\longrightarrow} \bigoplus_{x\in X} \Bbb Z,
\end{displaymath}
where the sum on the middle term ranges over all irreducible curves on $X$
and that on the right hand side over all closed points of $X$.
The map $\dt$ is the so-called tame symbol and $\dd$ is the sum of divisors
of rational functions on curves. We have the localization exact sequence
$$ CH^2(U,2) \to CH^1(Z,1) \to CH^2(X,1),$$
where
$$ CH^1(Z,1)=\Ker(\underset{1\leq i\leq 3}{\bigoplus} \Bbb C(Z_i)^*
\rmapo{\dd}  \underset{x\in Z}{\bigoplus} \Bbb Z)
\quad\text{ with } Z_i=X\cap Y_i.$$
By [AS2], Th.(6-1) we get the following.

\begin{thm}\label{0-5}
For $t\in M\setminus \MNL$, $\oXt$ is injective on
the subspace
$$\Sigma_t:= \Im(CH^1(\Zt,1) \to CH^2(\Xt,1))\otimes\Bbb Q\subset
CH^2(\Xt,1)\otimes\Bbb Q.$$
\end{thm}\medbreak

In \S 6 we show there exists $t\in M\setminus \MNL$ such that
$\Sigma_t\not=0$ so that Th.\ref{0-5} has a non-trivial implication on
the injectivity of $\oXt$. For this we need introduce some
special locus in the moduli space $M$.

\begin{defn}\label{d0-3}
Let $\Ta\subset M$ be the locus of those $X$ defined by an equation
$$ F=w A + z_1 z_2 B + c_1 z_1^d + c_2 z_2^d
\text{ for some } w\in P^1,\; A\in P^{d-1},\;B \in P^{d-2},\;
c_1,c_2 \in \Bbb C^*.$$
We define $\Tb$ (resp. $\Tc$) similarly by replacing
$(z_1,z_2)$ by $(z_2,z_3)$ (resp. $(z_3,z_1)$).
\end{defn}
\medbreak

We note that $\Tspqc\subset \Ta$ with $\sigma$, the identity, and $p=1,q=0$.
For $X$ in $\Ta$ defined by such an equation as above we consider
the following element
$$\caX=((\frac{z_2}{w})_{|Z_1},(\frac{w}{z_1})_{|Z_2},1)\in
\Bbb C(Z_1)^*\oplus \Bbb C(Z_2)^*\oplus \Bbb C(Z_3)^*.$$
It is easy to check $\caX\in CH^1(Z,1)$.
For $X$ in $\Tb$ (resp. $\Tc$) we define
an element $\cbX$ (resp. $\ccX$) in $CH^1(Z,1)$ by the same say.
Let $[\cijX]\in CH^2(X,1)$ be the
image of $\cijX\in CH^1(Z,1)$ for $(i,j)=(1,2)$ or $(2,3)$ or $(3,1)$.

\begin{thm}\label{0-6}
\begin{itemize}
\item[$(1)$]
If $d\geq 4$, $\Ta\not\subset \MNL$ and
$\oXt([\caXt])\not=0$ for $\forall t\in \Ta\setminus \MNL$.
\item[$(2)$]
If $d\geq 6$, $\Ta\cap\Tb\not\subset \MNL$ and
$\oXt([\caXt]),\; \oXt([\cbXt])$ are linearly independent
for $\forall t\in (\Ta \cap \Tb)\setminus \MNL$.
\item[$(3)$]
If $d\geq 10$, $\Ta\cap\Tb\cap\Tc\not\subset \MNL$ and
$\oXt([\caXt]),\; \oXt([\cbXt]),\; \oXt([\cbXt])$ are linearly independent
for $\forall t\in (\Ta \cap \Tb \cap \Tc)\setminus \MNL$.
\end{itemize}
\end{thm}\medbreak

The authors are grateful to Prof. A. Otwinowska for stimulating discussions
and her teaching us on her work [Ot].
\vskip 20pt

\section{Component of $\MNL$}

\medbreak

Let the notaion be as in Th.\ref{0-4}.
In what follows we fix $0\in M$ and let $X$ be the fibers over $0\in M$
of the universal family $\cX/M$ and write
$$ U=X\setminus Z,\; Z=\underset{1\leq j\leq 3}{\cup} Z_j\text{ with }
Z_j=X\cap Y_j.$$

\begin{defn}\label{d1-1} We put
$$ \aU =\{\frac{z_2}{z_1},\frac{z_3}{z_1}\}\in CH^2(U,2),$$
$$ \wU =\dlog{\frac{z_2}{z_1}} \wedge \dlog{\frac{z_3}{z_1}}
\in H^0(X,\WXZ 2),$$
where $z_j/z_i$ is viewed as an element of
$CH^1(U,1)=\Gamma(U,\cO^*_{U_{Zar}})$.
Note that $\wU = \regU 2(\aU)$, under the identification
$ F^2H^2(U,\Bbb C)\cong H^0(X,\WXZ 2)$.
\end{defn}

\begin{lem}\label{l1-1}
\begin{itemize}
\item[$(1)$]
$CH^2(U,2)_{dec}$ is generated by $\aU$ and $\{c,z_j/z_i\}$
with $c\in \Bbb C$ and $1\leq i,j\leq 3$.
\item[$(2)$]
$H^2(U,\Bbb Q(2))_{triv}=\Bbb Q\cdot \wU$.
\end{itemize}
\end{lem}

\proof The first assertion follows from [AS2] Lem.(2-1).
The second assertion follows from the fact
$H^2(U,\Bbb Q(2))_{triv}=\regU 2(CH^2(U,2)_{dec})$ by loc.cite.
\medbreak

Let $\Tspqc\subset M$ be as in Def.\ref{d0-2}.

\begin{lem}\label{l1-2}$\Tspqc$ is smooth irreducible and
$\codim_M(\Tspqc)=\binom{d+2}{2}-5$.
\end{lem}

\proof
Left to the readers as an easy exercise.
\medbreak

Assume $0\in \Tspqc$ and that $X$ is defined by such an equation as in
Def.\ref{d0-2}:
$$ F=w A +\prod_{1\leq \nu\leq r}
(c z_{\sigma(1)}^{p+q}-c_\nu z_{\sigma(2)}^p z_{\sigma(3)}^q).$$
We note that $w\not\in\sum_{1\leq j\leq 3} \Bbb C\cdot z_j$ by the
assumption that $X$ transversally intersects $Y$.

\begin{defn}\label{d1-2} We define
$$ \SU:=\Bbb C \cdot \wU \oplus \Bbb C \cdot \xU \subset
H^0(X,\WXZ 2)=F^0H^2(U,\Bbb C),$$
$$ \xU=\dlog{\frac{z_{\sigma(2)}^p z_{\sigma(3)}^q}{z_{\sigma(1)}^{p+q}}}
\wedge \dlog{\frac{w}{z_{\sigma(1)}}} \in H^0(X,\WXZ 2).$$
\end{defn}
\medbreak
We note that $\xU$ is apparently holomorphic only on $U\setminus W$,
where $W=U \cap \{w=0\}$ while it is easy to see that
its residue along any irreducible component of $W$ is zero.
Rewriting the equation of $X$ as
\begin{equation}\label{e1-1}
 w A + \prod_{1\leq \nu\leq r}
(c z_{\sigma(1)}^{p+q}-c_\nu z_{\sigma(2)}^p z_{\sigma(3)}^q)
= w A + \prod_{\mu\in I}
(c z_{\sigma(1)}^{p+q}-c_\mu z_{\sigma(2)}^p z_{\sigma(3)}^q)^{e_\mu},
\; (e_\mu\geq 1)
\end{equation}
where $c_\mu\not=c_{\mu'}$ if $\mu\not=\mu'\in I$,
$W$ is the disjoint sum of the following smooth irreducible components
for $\mu\in I$;
$$ \Wm=U\cap
\{w=c z_{\sigma(1)}^{p+q}-c_\mu z_{\sigma(2)}^p z_{\sigma(3)}^q=0\}.$$
We consider the condition:
\begin{equation}\label{i1-2}
\text{$\uc=[c'_\nu]_{1\leq \nu\leq r}$ such that
$c'_\nu$ is a root of unity for $1\leq \forall \nu \leq r$.}
\end{equation}

\begin{prop}\label{p1-1}
\begin{itemize}
\item[$(1)$]
  $ \SU\cap H^2(U,\Bbb Q(2)) =
\left.\left\{\begin{gathered}
   \Bbb Q \cdot \wU \oplus \Bbb Q \cdot \xU \\
   \Bbb Q \cdot \wU \\
   \end{gathered}\right.\qquad
\begin{aligned}
  &\text{if $\eqref{i1-2}$ holds}\\
  &\text{otherwise}\\
\end{aligned}\right.  $
\item[$(2)$]
$\SU \cap H^2(U,\Bbb Q(2)) \subset \Im(\regU 2)$.
\end{itemize}
\end{prop}

\begin{cor}\label{c1-2}
If the condition \eqref{i1-2} holds, $\Tspqc\subset \MNL$.
\end{cor}
\medbreak
The corollary follows immediately from Pr.\ref{p1-1} and Lem.\ref{l1-1}.
\medbreak

Noting $\wU \in H^2(U,\Bbb Q(2))$, Pr.\ref{p1-1} follows from
the following two claims.

\Claim {1} \it $\SU\cap H^2(U,\Bbb Q(2))\subset
\Bbb Q\cdot \wU \oplus \Bbb Q\cdot \xU$.
\rm

\Claim {2} \it $\xU \in H^2(U,\Bbb Q(2))$ if and only if \eqref{i1-2} holds,
in which case $\xU\in \Im(\regU 2)$.
\rm
\medbreak

We prove Claim 1.
For simplicity we assume that $\sigma\in \PG 3$ is the identity.
The following argument works in general case as well. Define
$$  Z_{ij}=X\cap \{z_i=z_j=0\}\quad (1\leq i\not=j\leq 3),$$
$$  Z_i=X\cap \{z_i=0\},\;
  V_i = Z_i \cap (\underset{1\leq j\not=i \leq 3}{\cup} Z_j).$$
We consider the composite map of the successive residue maps
$$\delta_{ij}: H^0(X,\WXZ 2) \rmapo{Res_{Z_i}}
H^0(Z_i,\Omega_{Z_i}(\log V_i)) \rmapo{Res_{Z_{ij}}}
H^0(Z_{ij},\cO_{Z_{ij}}) \cong \Psi\otimes\Bbb C,$$
where $\Psi=H^0(Z_{ij},\Bbb Q)=\bigoplus_{x\in Z_{ij}} \Bbb Q $.
For $\phi\in  H^0(X,\WXZ 2)$, $\phi\in H^2(U,\Bbb Q(2))$ implies
$\delta_{ij}(\phi) \in \Psi \subset \Psi\otimes\Bbb C$.
Now an easy residue calculation shows
$$ \delta_{12}(\wU)=\delta_{23}(\wU)=\delta_{31}(\wU)=- \underline{u},$$
$$ \delta_{12}(\xU)=-p \cdot\underline{u},\;
\delta_{31}(\xU)=q \cdot\underline{u},\;
\delta_{23}(\xU)=0,$$
where $\underline{u}=(1,1,\dots,1)\in \Psi$.
Thus, if $\phi=a \wU + b \xU\in H^2(U,\Bbb Q(2))$ with $a,b\in \Bbb C$,
it implies
$-(a+bp),-a,-(a-bq)\in \Bbb Q$.
Noting that at least one of $p$ and $q$ is not zero, it implies $a,b\in \Bbb Q$
and the proof of Claim 1 is complete.
\medbreak

Next we prove Claim 2. Consider
$$ \beta=\{\frac{z_{\sigma(2)}^p z_{\sigma(3)}^q}{z_{\sigma(1)}^{p+q}},
\frac{w}{c z_{\sigma(1)}}\} \in CH^2(U',2)\quad (U':=U\setminus W)$$
We have the commutative diagram
$$
\begin{matrix}
CH^2(U,2)\otimes \Bbb Q & \rmapo{\regU 2} &
F^0H^2(U,\Bbb Q(2)) & \rmapo{\hookrightarrow} & F^2H^2(U,\Bbb C) \\
\downarrow && \downarrow\rlap{$\iota_1$} && \downarrow\rlap{$\iota_2$} \\
CH^2(U',2)\otimes \Bbb Q & \rmapo{\regUd 2} &
F^0H^2(U',\Bbb Q(2)) & \rmapo{\hookrightarrow} & F^2H^2(U',\Bbb C)
\end{matrix}
$$
and we have
$\regUd 2(\beta)=\iota_2(\xU)$ in $F^2H^2(U',\Bbb C)$. Since $\iota_2$ is
injective, the first part of the claim follows from the following assertion:
\begin{equation}\label{i1-3}
\text{
$\regUd 2(\beta)\in \Im(\iota_1)$ if and only if \eqref{i1-2} holds.}
\end{equation}
To show this we consider the commutative diagram
$$
\begin{matrix}
  CH^2(U,2)\otimes \Bbb Q &\longrightarrow&
CH^2(U',2)\otimes \Bbb Q  &\rmapo{\partial_1}& CH^1(W,1)\otimes \Bbb Q \\
\downarrow\rlap{$\regDU 2$} && \downarrow\rlap{$\regDUd 2$} &&
\downarrow\rlap{$\regDW 1$}\\
 H^2_D(U,\Bbb Q(2)) &\longrightarrow& H^2_D(U',\Bbb Q(2)) &\longrightarrow&
 H^1_D(W,\Bbb Q(1))
\end{matrix}
$$
where $reg_{D,*}^*$ denotes the regulator map to Deligne cohomology.
We have the commutative diagram (cf. [EV])
$$
\begin{matrix}
0 \to H^1(U,\Bbb Q(1))\otimes\Bbb C/\Bbb Q(1) &\longrightarrow&
 H^2_D(U,\Bbb Q(2)) &\rmapo{\pi_U}& F^0H^2(U,\Bbb Q(2)) \to 0 \\
\downarrow && \downarrow && \downarrow\rlap{$\iota_1$} \\
0 \to H^1(U',\Bbb Q(1))\otimes\Bbb C/\Bbb Q(1) &\longrightarrow&
 H^2_D(U',\Bbb Q(2)) &\rmapo{\pi_{U'}}& F^0H^2(U',\Bbb Q(2)) \to 0 \\
 \downarrow\rlap{$\psi$} &&  \downarrow &&  \downarrow\rlap{$\partial_2$}\\
0 \to H^0(W,\Bbb Q)\otimes\Bbb C/\Bbb Q(1)  &\longrightarrow&
 H^1_D(W,\Bbb Q(1))   &\rmapo{\pi_W}&  F^0H^1(W,\Bbb Q(1)) \to 0
\end{matrix}
$$
where the composite of $reg_{D,*}^*$ with $\pi_*$ conincide with the regulator
map to singular cohomology. The horizontal sequences are exact. The vertical
sequences are localization sequences and they are exact except the most right
one. We have (cf. \eqref{e1-1})
$$ \Coker(\psi)\cong \Bbb C/\Bbb Q(1)\otimes\Phi, \quad
\Phi:=\Coker(\Bbb Q \to \bigoplus_{\mu\in I}
\Bbb Q;\; 1 \to (e_\mu)_{\mu\in I}).$$
To see this we note the commutative diagram
$$
\begin{matrix}
CH^1(U',1)\otimes\Bbb Q &\rmapo{\phi}& \bigoplus_{\mu\in I} \Bbb Q \\
\downarrow\rlap{$\regUd 1$} && \downarrow\rlap{$\cong$}\\
H^1(U',\Bbb Q(1)) &\longrightarrow& H^0(W,\Bbb Q)
\end{matrix}
$$
where $\phi$ is given by taking orders of fuctions along the compoenents of
$W$. One easily sees that $\regUd 1$ is surjective and that
$CH^1(U',1)$ is generated by $\Bbb C^*$, $w/z_1$ and $z_j/z_i$ with
$1\leq i,j\leq 3$ and the desired assertion follows.
Thus the above diagram gives rise to the exact sequence
$$ F^0H^2(U,\Bbb Q(2)) \overset{\iota_1}{\longrightarrow}
F^0H^2(U',\Bbb Q(2))_{\partial}
\overset{\delta}{\longrightarrow}  \Bbb C/\Bbb Q(1)\otimes\Phi,$$
where
$F^0H^2(U',\Bbb Q(2))_{\partial}=\Ker(\partial_2)$.
Now an easy calculation shows
$\partial_1(\beta)=(c_{\mu}^{-e_\mu})_{\mu\in I}$.
Noting the commutative diagram
$$
\begin{matrix}
0 \to & \Bbb C^*\otimes\Bbb Q  &\longrightarrow& CH^1(\Wm,1)\otimes\Bbb Q \\
&\downarrow\rlap{$\log$} && \downarrow\rlap{$\regDWm 1$} &&\\
0 \to & \Bbb C/\Bbb Q(1) &\longrightarrow& H_D^1(\Wm,\Bbb Q(1))
&\rmapo{\pi_{\Wm}}& F^0H^1(\Wm,\Bbb Q(1))\to 0
\end{matrix}
$$
it implies $\regUd 2(\beta) \in  F^0H^2(U',\Bbb Q(2))_{\partial}$ and
$$ \delta(\regUd 2(\beta))=
\text{ the class of $(-e^{\mu}\log{c_\mu})_{\mu\in I}$ in $\Phi$}.$$
This proves \eqref{i1-3}.
In order to show the second part of Claim 2, assume
that there is $c'\in \Bbb C^*$ such that $c'_\nu:=c c_\nu $ is a root of
unity for $1\leq \forall \nu\leq r$. Taking
$$ \beta'=\{\frac{z_{\sigma(2)}^p z_{\sigma(3)}^q}{cc' z_{\sigma(1)}^{p+q}},
\frac{w}{z_{\sigma(1)}}\} \in CH^2(U',2),$$
$\regUd 2 (\beta') =\regUd 2(\beta)=\iota_2(\xU)$ and
$\partial_1(\beta)=({c'_{\mu}}^{-e_\mu})_{\mu\in I}=0\in
CH^1(W,1)\otimes\Bbb Q$.
It implies that $\beta'$ has a lift $\beta''\in CH^2(U,2)\otimes\Bbb Q$.
Then $\regU 2(\beta'')=\xU\in F^0H^2(U,\Bbb Q(2))$
by the injectivity of $\iota_2$.
This completes the proof of Claim 2.
\vskip 20pt

\section{Infinitesimal interpretation}

\medbreak

In this section we take the first step of the proof of Th.\ref{0-4}.
Let the assumption and the notation be as in \S1.
Take $\Delta\subset M$, a simply connected neighbouthood of $0$ in $M$.
For $\lambda\in H^2(U,\Bbb C)$ and $t\in \Delta$,
let $\lt\in H^2(\Ut,\Bbb C)$ be the flat translation of $\lambda$
with respect to the Gauss-Maninn connection
$$ \nabla : \HOcU 2 \to \WM 1 \otimes \HOcU 2,$$
where $\HOcU p$ is the sheaf of holomorphic sections of the local system
$\HCcU p :=R^p f_* \Bbb C$ with $f:\cU \to M$, the natural morphism.
We sometime consider $\lambda$ a section over $\Delta$ of $\HCcU  2$ via
$H^2(U,\Bbb C)\cong \Gamma(\Delta,\HCcU 2)$. Put
$$ \Dl=\{t\in \Delta|\; \lt \in F^2H^2(\Ut,\Bbb C)\}.$$
$\Dl$ is a closed analytic subset of $\Delta$ since it is defined by the
vanishing of the image of $\lambda$ under the map
$$ \Gamma(\Delta,\HCcU 2) \to \Gamma(\Delta,\HOcU 2/F^2\HOcU 2)$$
where $F^q\HOcU p\subset \HOcU p$ is the Hodge subbundle.
Taking $\Delta$ sufficiently small if necessary,
we have by Lem.\ref{l1-1}(2)
\begin{equation}
\label{i2-1}
\text{$ \MNL \cap \Delta = {\cup}_{\lambda} \Dl$ with
$\lambda$ ranging over $F^0H^2(U,\Bbb Q(2)) \setminus \Bbb Q\cdot \wU$.}
\end{equation}
By Griffiths transversality, $\nabla$ induces
$$  \nabb \scs \HU 2 0 \to \WMo {1}\otimes \HU {1}{1},$$
where for integers $p,q$ we put
$\HU p q=H^q(X,\WXZ p)$ and $\WMo 1$ is the fiber of $\WM 1$ at $0\in M$.
Let $\TMo$ be the tangent space of $M$ at $0$.
Via the natural isomorphism $\TMo\cong \Hom(\WMo 1,\Bbb C)$, it induces
the pairing
$$ <\;,\;> \scs \TMo \otimes \HU 20 \to \HU 11.$$
For $\lambda\in \HU 20$ write
$$ \Vl:=\{\partial \in \TMo|\; <\partial,\lambda>=0\}.$$
By the construction we have
\begin{equation}
\label{i2-2}
T_0(\Dl) \subset \Vl.
\end{equation}

\begin{thm}\label{t2-1}Let $\lambda\in \HU 20$ and assume
$\lambda \not\in \Bbb C\cdot \wU$.
\begin{itemize}
\item[$(1)$]
$\codim_{\TMo}(\Vl)\geq \binom{d+2}{2}-5$.
\item[$(2)$]
Assume $d\geq 4$ and that $T\subset \Dl$ is an ireducible component of
$\codim_{\Delta}(T) = \binom{d+2}{2} -5$.
Then $T=\Dl$, and $0\in \Tspqc\cap \Delta$
for some $\sigma,p,q,\uc$, and $\lambda\in \SU$.
Moreover, if $\lambda\in H^2(U,\Bbb Q(2))$, then
$\Dl= \Tspqc\cap \Delta$.
\end{itemize}
\end{thm}
\medbreak

\begin{thm}\label{t2-2}
Assume $0\in \Tspqc$ for some $\sigma,p,q,\uc$. Then we have
$$ \SU =\{\omega\in \HU 20|\; <\partial,\omega>=0\text{ for }
\forall \partial\in T_0(\Tspqc) \}$$
\end{thm}
\medbreak

In the rest of this section we deduce Th.\ref{0-4}
from Th.\ref{t2-1} and Th.\ref{t2-2}.
Th.\ref{0-4}(1) and (2)
follow immediately from Th.\ref{t2-1} and \eqref{i2-1} and Pr.\ref{p1-1}
by noting
$$ \codim_{\Delta}(\Dl)\geq\codim_{T_0(\Delta)}(T_0(\Dl))
\geq\codim_{T_0(\Delta)}(\Vl)$$
where the last inequality is due to \eqref{i2-2}.
Assume $0\in \Tspqc$. We shall show that there exists a subset
$E\subset \Delta_T:=\Tspqc\cap \Delta$ which is the union of a countable
number of proper closed analytic subsets of $\Delta_T$ such that
$F^0H^2(\Ut,\Bbb Q(2))\subset \SUt$ for $\forall t\in \Delta_T \setminus E$.
By Pr.\ref{p1-1} the last condition implies
$F^2H^2(\Ut,\Bbb Q(2))\subset \Im(\regUt 2)$. Hence Th.\ref{0-4}(3) follows.
Write
$H^2(U,\Bbb Q(2))=\{\lambda_i\}_{i\in I}$ as a set and put
$$A=\{i\in I|\; \Delta_T \subset \Dli \},\;
B=\{i\in I|\; \Delta_T \not\subset \Dli\},\;
E=\Delta_T\cap (\underset{i\in B}{\cup}  \Dli).$$
Note that $I$ is countable and $I=A\cup B$ and $A\cap B=\emptyset$.
For $\forall t\in \Delta_T-E$, we have
$F^0H^2(\Ut,\Bbb Q(2))=\{\lit\}_{i\in A}$ so that
$H^2(\Ut,\Bbb C)\isom \Gamma(\Delta_T,\HCcU 2)$ induces
$$
   F^2H^2(\Ut,\Bbb Q(2))\hookrightarrow
   \Gamma(\Delta_T,\HCcU 2 \cap F^2\HOcU 2),
$$
which further implies
$$ F^2H^2(\Ut,\Bbb Q(2))\subset
\Ker(\HUt 20 \to \Omega^1_{\Delta_T,t}\otimes\HUt 11).$$
Th.\ref{t2-2} implies that the last space is equal to $\SUt$
and the desired assertion follows. This completes the proof of Th.\ref{0-4}.
\vskip 20pt

\section{Reduction to Jacobian rings}

\medbreak

Let the assumption be as in \S2. In this section we rephrase the theorems
in \S2 in terms of Jaocibian rings and
prove Th.\ref{t2-1}(1) and Th.\ref{t2-2}.
Let $P=\Bbb C[z_0,z_1,z_2,z_3]$ be the homogeneuous coordinate ring of
$\Bbb P^3$. For an integer $l>0$ let $P^l\subset P$
be the subspace of homogeneous polynomials of degree $l$.
Let the assumption be as in \S2 and fix $F\in P^d$ which defines
$X\subset \Bbb P^3$. Consider the ideal
of $P$;
$$ J_F=<\frac{\partial F}{\partial z_0},
z_1 \frac{\partial F}{\partial z_1}, z_2 \frac{\partial F}{\partial z_2},
z_3 \frac{\partial F}{\partial z_3}>.$$
The assumption that $X$ transversally intersects $Y$ is equivalent to
the condition:
\begin{itemize}
\item[\bf(3-1)]
$J_F$ is complete intersection of degree $(d-1,d,d,d)$.
\end{itemize}
Write
$$ R_F=P/J_F, \; \JF l =J_F\cap P^l,\;
\RF l =\Im(P^l \to R_F)=P^l/\JF l.$$
Note $F\in \JF d$. We have the following well-known facts:
\begin{enumerate}
\item[\bf(3-2)]
We have the canonical surjective homomorphism
 $$ \psi \scs P^d \longrightarrow \TMo\; ; \;
G \to \{F+\epsilon G = 0\}\subset\Bbb P^3_{\Bbb C[\epsilon]},$$
where $\Bbb C[\epsilon]$ is the ring of dual numbers.
We have $\Ker(\psi)=\Bbb C\cdot F$.
\item[\bf(3-3)]\label{i3-3}
We have the isomophisms
$$ \phi: P^{d-1} \isom \HU 20,\quad
\phi': \RF {2d-1} \isom \HU 11,$$
such that the diagram
  $$
    \begin{matrix}
      P^d\otimes P^{d-1}  &\overset{\mu}{\longrightarrow}&\RF {2d-1} \\
      \downarrow\rlap{$\psi\otimes\phi$} &&
      \downarrow\rlap{$\phi'$} \\
      \TMo \otimes \HU 20  & \overset{<\;,\;>}{\longrightarrow} & \HU 11 \\
   \end{matrix}
  $$
commutes up to non-zero scalar where $\mu$ is the multiplication.
\item[\bf(3-4)]
We have the following formula
$$ \phi(G) = Res_X \frac{G}{z_1z_2z_3} \Omega \quad (G\in P^{d-1}),$$
where
$\Omega=\sum_{i=0}^{3} (-1)^i z_i dz_0\wedge \cdots \wedge
\widehat{dz_i}\wedge \cdots \wedge dz_3\in H^0(\Bbb P^3,\Omega^3_{\Bbb P^3})$
and
$$ Res_X \scs H^0(\Bbb P^3,\Omega^3_{\Bbb P^3}(\log Y)) \to
H^0(X,\WXZ 2)\quad (Y=\{z_1z_2z_3=0\}\subset \Bbb P^3)$$
is the residue map.
\end{enumerate}
We omit the proof of the following lemma which can be easily shown by using
(3-4).

\begin{lem}\label{l3-1}
\begin{itemize}
\item[$(1)$] Putting $\wF=\frac{\partial F}{\partial z_0}$,
$\phi(\wF)=\wU$ (cf. Def.\ref{d1-1}).
\item[$(2)$]
Assume $0\in \Tspqc$ and that $X$ is defined by such an equation
as Def.\ref{d0-2}:
$$ F=w A +\prod_{1\leq \nu\leq r}
(c z_{\sigma(1)}^{p+q}-c_\nu z_{\sigma(2)}^p z_{\sigma(3)}^q).$$
Put
$$ \xF =
\pd{w}{z_0} (q z_{\sigma(2)}\frac{\partial A}{\partial z_{\sigma(2)}} -
p z_{\sigma(3)}\frac{\partial A}{\partial z_{\sigma(3)}}) -
\frac{\partial A}{\partial z_0}(q z_{\sigma(2)} \pd{w}{z_{\sigma(2)}} -
p z_{\sigma(3)} \pd{w}{z_{\sigma(3)}})\; \in P^{d-1}.$$
Then we have
$\phi(\xF)=\xU$ (cf. Def.\ref{d1-2}).
\end{itemize}
\end{lem}
\medbreak

In what follows we identify $P^{d-1}=\HU 20$ via $\phi$.
For $\lambda \in P^{d-1}$ write
$$ \Il d=\{x\in P^d|\; \lambda x=0 \in \RF {2d-1}\}.$$
\eqref{i2-2} implies
\begin{description}
\item[(3-5)]
$ \psi^{-1}(T_0(\Dl))\subset \Il d.$
\end{description}

In view of the above lemmas Th.\ref{t2-1}(1) and Th.\ref{t2-2}
follow from the
following theorems.

\begin{thm}\label{t3-1}Assuming $\lambda\not\in \JF {d-1}=\Bbb C\cdot \wF$,
$\dim(P^d/\Il d)\geq \binom{d+2}{2} -5$. The equality holds if and only if
$\Il d$ is complete intersection of degree $(1,d-1,d,d)$.
\end{thm}

\begin{thm}\label{t3-2}Let the assumption be as in Lem.\ref{l3-1}(2).
\begin{itemize}
\item[$(1)$]
$\psi^{-1}(T_0(\Tspqc)) = wP^{d-1} + \JF d.$
\item[$(2)$]
$\Il d= wP^{d-1} + \JF d $
if $\lambda=a\wF + b\xF$ with $b\not=0$.
\item[$(3)$]
$\Bbb C\cdot \wF \oplus \Bbb C\cdot \xF =
\{\lambda \in P^{d-1}|\; \lambda x=0 \in \RF {2d-1}
\text{ for } \forall x\in wP^{d-1} + \JF d \}.$
\end{itemize}
\end{thm}

\medbreak

In the rest of this section we prove Th.\ref{t3-1} and Th.\ref{t3-2}.
We need the following theorems. The first one is Macaulay's theorem and
we refer [GH], p659, for the proof.
The second one is due to A. Otwinowska and is shown
by the same method as the proof of [Ot], Th.2.

\begin{thm}\label{t3-3}
There exists a natural isomorphism
$$ \tau_F \scs \RF {4d-5} \isom \Bbb C$$
and the pairing induced by multiplication
$$ \RF l \otimes \RF {4d-5-l} \to \RF {4d-5}
\overset{\tau_F}{\longrightarrow} \Bbb C$$
is perfect for $\forall l$.
\end{thm}\medbreak

\begin{thm}\label{t3-4}
Let $I\subset P$ be a homogeneous ideal satisfying
the conditions:
\begin{itemize}
\item[$(1)$]
$I$ is Gorenstein of degree $N>0$, namely there exists a non-zero
linear map $\mu: P^N \to \Bbb C$ such that
$I^l=\{x\in P^l|\; \mu(xy)=0 \text{ for } \forall y\in P^{N-l}\}.$
\item[$(2)$]
$I$ contains a homogeneous ideal $J$ which is complete intersection of degree
$(e_0,e_1,e_2,e_3)$ with $e_0\leq e_1 \leq e_2 \leq e_3$.
\item[$(3)$]
There is an integer $b$ such that $e_0\leq b\leq e_1-1$ and
$N+3=e_2+e_3+b$.
\end{itemize}
For $\forall l\geq 1$ we have
$$ \dim(P^l/I^l)\geq \dim(P^l/<z_0,z_1^b,z_2^{e_2},z_3^{e_3}>\cap P^l).$$
Moreover, if the equality holds for some $l_0\leq N-b$,
then there is a complete intersection ideal $I_0$
of degree $(1,b,e_2,e_3)$ such that $I^l=I^l_0$ for all $l\leq l_0$.
\end{thm}\medbreak

Now we start the proof of Th.\ref{t3-1}.
For $\lambda \in P^{d-1}$ consider the linear map
$$ \lambda^* \scs P^{3d-4} \to \Bbb C\; ; \; x \to \tau_F(\lambda x).$$
For an integer $l\geq 0$ define
$$\Il l=\{x\in P^l|\; \lambda^*(xy)=0 \text{ for } \forall y\in P^{3d-4-l}\}.$$
By Th.\ref{t3-3} $\lambda^*\not=0$
if and only if $\lambda\not\in \JF {d-1}$ and
$\Il l$ in case $l=d$ coincides with $\Il d$ defined before.
Define a homogeneous ideal of $P$ by
$$ \Ill =\underset{l\geq 0}{\bigoplus} \Il l \subset P.$$
We take $N=3d-4$, $(e_0,e_1,e_2,e_3,b)=(d-1,d,d,d,d-1)$
and apply Th.\ref{t3-3} to
$I=\Ill$ and $J=J_F$ noting (3-1). Since
$$\dim(P^d/<z_0,z_1^{d-1},z_2^{d},z_3^{d}>\cap P^d) =\binom{d+2}{2}-5,$$
it implies Th.\ref{t3-1}.
\medbreak

Next we show Th.\ref{t3-2}.
Let $\PGL_4$ be the group of projective transformations on $\Bbb P^3$
and let $G\subset \PGL_4$ be the subgroup of such $g\in \PGL_4$ that
$g(Y_j)=Y_J$ for $1\leq \forall j\leq 3$.
It is evident that $G$ naturally acts on $M$ and $\Tspqc\subset M$ is stable
under the action. Let ${\Tspqc}_{(w,c)}\subset \Tspqc$ be the closed subset of
those surfaces defined by equations of the form
$$ w B +\prod_{1\leq \nu\leq r}
(c z_{\sigma(1)}^{p+q}-c_\nu z_{\sigma(2)}^p z_{\sigma(3)}^q)
\quad\hbox{ for some } B\in P^{d-1}.$$
It is easy to see that the natural map
$G\times {\Tspqc}_{(w,c)} \to \Tspqc$ is smooth and surjective and that
$ \psi^{-1}(T_0({\Tspqc}_{(w,c)}))=w P^{d-1}$.
The map
$\TMo \rmapo{\psi^{-1}} P^d/\Bbb C\cdot F \to \RF d$
identifies $\RF d$ with the quotient of $\TMo$ by the infinitesimal action of
the tangent space at the identity of $G$. It implies
$$ \psi^{-1}(T_0(\Tspqc))=\pi^{-1}\pi\psi^{-1}(T_0({\Tspqc)}_{(w,c)}))=
w P^{d-1} +\JF d.$$
This completes the proof of Th.\ref{t3-2}(1).
\medbreak

For $\lambda=a\wF + b\xF$, an easy calculation shows
$\lambda w \in \JF d$ so that
$\Il d\supset wP^{d-1} +\JF d$. Assuming $b\not=0$ we have
\begin{multline*}
\binom{d+2}{2}-5 \underset{(*)}{\leq} \dim(P^d/\Il d) \leq
\dim(P^d/wP^{d-1} +\JF d) \\
\underset{(**)}{=} \codim_{\TMo}(T_0(\Tspqc))
\leq \codim_{M}(\Tspqc) \underset{(***)}{\leq} \binom{d+2}{2}-5,
\end{multline*}
where $(*)$ follows from Th.\ref{t3-1}, $(**)$ from Th.\ref{t3-2}(1), and
$(***)$ from Lem.\ref{l1-2}.
Thus the above inequalities are all equalities so that
$\Il d =wP^{d-1} +\JF d$. This completes the proof of Th.\ref{t3-2}(2).
\medbreak

For $\lambda \in P^{d-1}$ we have the following equivalences
\begin{align*}
\lambda w y =0\in \RF {2d-1} \text{ for } \forall y\in P^{d-1}
& \Leftrightarrow \lambda w y z =0\in \RF {4d-5} \text{ for }
\forall y\in P^{d-1}, \; \forall z\in P^{2d-4} \\
& \Leftrightarrow \lambda w x =0\in \RF {4d-5} \text{ for }
\forall x\in P^{3d-5} \\
& \Leftrightarrow \lambda w =0\in \RF {d} \\
\end{align*}
where the first and the last euivalences follows from Th.\ref{t3-3}.
Hence it suffices to show
$$ \Bbb C\cdot \wF \oplus  \Bbb C\cdot \xF =
\Ker(P^{d-1} \overset{w}{\longrightarrow} \RF d).$$
We have already senn that the left hand side is contained in the righ hand
side. Thus it suffices to show
$\dim(\Ker(P^{d-1} \overset{w}{\longrightarrow} \RF d)=2$.
We have
$\dim(P^{d-1})=\binom{d-1+3}{3}$.
By (3-1) we have
$$ \dim(\RF d)=\dim(P^d/P^d\cap <z_0^{d-1},z_1^d,z_2^d,z_3^d>)=
\binom{d+3}{3}-
\left(
\binom{1+3}{3}+3\binom{0+3}{3}
\right).$$
We easily see that $<w>+J_F$ is complete intersection of degree
$(1,d-1,d,d)$ so that
\begin{multline*}
\dim(\Coker(P^{d-1} \overset{w}{\longrightarrow} \RF d))
 = \dim(P^d/P^d\cap <z_0,z_1^{d-1},z_2^d,z_3^d>) \\
 = \binom{d+3}{3}-\left(
 \binom{d-1+3}{3}+\binom{1+3}{3}+2\binom{0+3}{3}\right)
+\binom{0+3}{3}
\end{multline*}
These imply the desired assertion and the proof of Th.\ref{t3-2}
is complete.
\vskip 20pt

\section{Proof of key theorem}

\medbreak

In this and next sections we prove Th.\ref{t2-1}(2) to complete the proof of
Th.\ref{0-4}. Let the assumption be as in Th.\ref{t2-1}(2).
For $t\in \Delta$ let $\Ft\in P^d$ define
$\Xt\subset \Bbb P^3$, $R_{\Ft}$ be the corresponding Jacobian ring.
For $t\in \Dl$ let $I_{\lt}\subset P$ be defined in the same manner
as $\Ill$ with $\lambda$ replaced by $\lt\in \HUt 20$, the flat translation of
$\lambda$. For $\forall t\in T$ we have
$$ \codim_{\Delta}(T)\geq \codim_{T_t(\Delta)}(T_t(\Dl)) \geq
\dim (P^d/\Ilt d)\geq \binom{d+2}{2}-5,$$
where $T_t(*)$ denotes the tangent space at $t$.
The second inequality follows from (3-5) and the last from Th.\ref{t3-1}.
Hence
the assumption implies that the above inequalities are all equalities,
which implies $T=\Dl$ and $\psi^{-1}(T_t(\Dl))=\Ilt d$.
It also implies that $I_{\lt}$ is complete intersection
of degree $(1,d-1,d,d)$ so that $\Ilt 1=\Bbb C\cdot w_t$ for some
$w_t\in P^1$ determined up to non-zero scalar. It gives rise to the morphism
\begin{equation}\label{(4-1)}
h \scs \Dl \to \dP(P^1)=\Bbb P^3\; ; \;
t \to [w_t]:=\Bbb C\cdot w_t.
\end{equation}
For $[w]\in \dP(P^1)$ define the closed subset $\Dlw =h^{-1}(w)\subset \Dl$.
Note $ \codim_{\Dl}(\Dlw)\leq 3.$
In what follows we put $[w]=h(0) \in \dP(P^1)$.
For an integer $l\geq 0$ put $\Pw l=P/wP^{l-1}$. For $t\in M$ put
$$ \Phi_t =\Im(\JFt d \to \Pw d)=
\Im(\sum_{i=1}^3 \Bbb C\cdot z_i\frac{\partial \Ft}{\partial z_i}+
P^1\cdot \frac{\partial \Ft}{\partial z_0} \to \Pw d).$$
If $t\in \Dlw$, $I_{\lt}$ is complete intersection of degree $(1,d-1,d,d)$
with $\Ilt 1=\Bbb C\cdot w$. Noting $\JFt d\subset \Ilt d$, it implies
$$ \dim(\Phi_t)\leq \dim(\Ilt d/wP^{d-1})=\dim(\Pw 1)+2\dim(\Pw 0)=5.$$
Thus we have $\Dlw\subset M_w$ where
$M_w=\{t\in M|\; \dim(\Phi_t)\leq 5\}$, which is a closed algebraic subset of
$M$. Put $E_w=\Bbb C^3\oplus \Pw 1$. For
$\Gamma=[\gamma_1:\gamma_2:\gamma_3:\Lq] \in \dP(E_w)$ put
$$ \Mwq=\{t\in M|\;
\sum_{i=1}^3\gamma_i z_i\frac{\partial \Ft}{\partial z_i}+
\Lq \frac{\partial \Ft}{\partial z_0} \in wP^{d-1}\}.$$
We note
$ M_w =\cup_{\Gamma} \Mwq$
with $\Gamma$ ranging over $\dP(E_w)$.

\begin{lem}\label{l4-1}
\begin{itemize}
\item[$(1)$]
If $\Mwq\not=\emptyset$, $\codim_{M_w}(\Mwq)\leq 5$.
\item[$(2)$]
For $t \in \Mwq$,
$\psi^{-1}(T_t(\Mwq))=\{G\in P^d|\;
\sum_{i=1}^3\gamma_i z_i\frac{\partial G}{\partial z_i}+
\Lq \frac{\partial G}{\partial z_0} \in wP^{d-1}\}.$
\end{itemize}
\end{lem}

\proof
Lem.\ref{l4-1}(2) follows directly from (3-2).
We prove Lem.\ref{l4-1}(1).
For $t\in M$ put
$$ \Qt=\{(\gamma_1,\gamma_2,\gamma_3,L)\in E_w|\;
\sum_{i=1}^3\gamma_i z_i\frac{\partial \Ft}{\partial z_i}+
\Lq \frac{\partial \Ft}{\partial z_0} \in wP^{d-1}\}.$$
It is a linear subspace of $E_w$ and $\dim(\Qt)=6-\dim(\Phi_t)$.
For an integer $0\leq e\leq 5$ put
$ \Mwe=\{t\in M_w|\; \dim(\Phi_t)=5-e\}$
which is a locally closed subset of $M_w$.
Letting $\Ge$ be the Grassman variety of $(e+1)$-dimensional linear subspaces
in $E_w$, we define the morphism
$\pie: \Mwe \to \Ge$ by $\pie(t)=\Qt$.
For $\Gamma\in \Go =\dP(E_w)$ put
$ \Geq =\{\Gamma'\in \Ge|\; \Gamma'\supset \Gamma\}$.
By definition we have set-theoretically
$$ \Mwq = \underset{0\leq e\leq 5}{\coprod}
\Mwe \times_{\Ge} \Geq.$$
It proves Lem.\ref{l4-1}(1) by noting that
$\Geq$ is smooth of $\codim_{\Ge}(\Geq)=5-e$.
\medbreak

Now we fix $\Gamma=[\gamma_1:\gamma_2:\gamma_3:\Lq] \in \dP(E_w)$ such that
$0\in \Mwq$. It means
\begin{description}
\item[(4-2)]
$ \sum_{i=1}^3\gamma_i z_i\frac{\partial F}{\partial z_i}+
\Lq \frac{\partial F}{\partial z_0} \in wP^{d-1}.$
\end{description}
Put $\Dlwq=\Dlw\cap \Mwq$. By Lem.\ref{l4-1}(1) and the fact
$\codim_{\Dl}(\Dlw)\leq 3$ (cf. \eqref{(4-1)}), we have
$$ 8=3+5\geq \codim_{\Dl}(\Dlwq)\geq \codim_{T_0(\Dl)}(T_0(\Dlwq)).$$
It implies that
$T_0(\Mwq)$ contains a subspace of codimension$\leq 8$ in $T_0(\Dl)$.
Recall that we have shown $\psi^{-1}(T_0(\Dl))=\Il d \supset w P^{d-1}$.
Hence Lem.\ref{l4-1}(2) implies that there exists a subspace $Q$ of
codimension$\leq 8$ in $P^{d-1}$ such that
$$ G(\sum_{i=1}^3\gamma_i z_i  \frac{\partial w}{\partial z_i}+
\Lq \frac{\partial w}{\partial z_0}) \in wP^{d-1}
\text{ for } \forall G\in Q.$$
If $\sum_{i=1}^3\gamma_i z_i  \frac{\partial w}{\partial z_i}+
\Lq \frac{\partial w}{\partial z_0}\not\in \Bbb C\cdot w$, it implies
$G\in w P^{d-2}$. Since
$\codim_{P^{d-1}}(w P^{d-2})=\binom{d+1}{2}$, this is a contradiction if
$\binom{d+1}{2}>8$ which holds when $d\geq 4$. Thus we get the condition:
\begin{description}
\item[(4-3)]
$\sum_{i=1}^3\gamma_i z_i  \frac{\partial w}{\partial z_i}+
\Lq \frac{\partial w}{\partial z_0} \in \Bbb C\cdot w$.
\end{description}
Now a key lemma is the following.

\begin{lem}\label{l4-2}There exists $t\in \Dl$ such that
$w_t\not\in \sum_{i=1}^3 \Bbb C\cdot z_i$.
\end{lem}\medbreak

We will prove Lem.\ref{l4-2} in the next section.
Admitting Lem.\ref{l4-2}, we finish the proof of Th.\ref{t2-1}(2).
Let
$$ \Dlo=\{t\in \Dl|\; w_t\not\in \sum_{i=1}^3 \Bbb C\cdot z_i\}.$$
By Lem.\ref{l4-2} it is a non-empty open subset of $\Dl$.
We may assume $0\in \Dlo$.
By transforming by an element of $G$ (cf. the proof of Th.\ref{t3-1}),
we may
suppose $w=z_0$. The condition (4-3) now reads
$\Lq\in \Bbb C\cdot z_0$. Then $\gamma_1,\gamma_2,\gamma_3$ are
not all zero and the condition (4-2) implies
$$ \sum_{i=1}^3\gamma_i z_i\frac{\partial F}{\partial z_i} \in z_0 P^{d-1}.$$
Writing
$F=z_0 B + C$ with $C$, a homogeneous polynomial of degree $d$ in
$\Bbb C[z_1,z_2,z_3]$, the above condition is equivalent to
$$ \sum_{i=1}^3\gamma_i z_i\frac{\partial C}{\partial z_i} =0.$$
Write
$$ C=\underset{\alphab=(\alpha_1,\alpha_2,\alpha_3)}{\sum} c_{\alphab}
z^{\alphab},\quad (z^{\alphab}= z_1^{\alpha_1}z_2^{\alpha_1}z_2^{\alpha_2}, \;
c_{\alphab}\in \Bbb C)$$
and take $\alphab$ with $c_{\alphab}\not=0$. The above condition implies that
$\alphab$ is an integral point lying on the sectional line $\ell$ in
$x_1x_2x_3$-space defined by
$$\ell\scs \sum_{i=1}^3 x_i-d = \sum_{i=1}^3\gamma_i x_i=0,
\; x_i\geq 0\; (i=1,2,3).$$
Furthermore the condition (3-1)
implies that $C$ is divisible by neither of
$z_1,z_2,z_3$. Writing $\pi_i : x_i=0$,
it implies that $\ell$ and $\pi_i$ intersect at an integral point for
$1\leq \forall i\leq 3$. This implies that $\ell$ passes through one of the
points $(d,0,0)$, $(0,d,0)$, $(0,0,d)$.
Assuming that $\ell$ passes through the first point, we get $\gamma_1=0$
and hence
$\alpha_2:\alpha_3=-\gamma_3:\gamma_2=p:q$
for some coprime non-negative integer $p,q$.
Writing $\alpha_2=p j,\; \alpha_3=q j$ with $j\in \Bbb Z$, we get
$\alpha_1=d-(p+q)j$ since $\sum_{i=1}^3 \alpha_i=d$.
The condition that $\ell$ and $\pi_1$ intersect at an integral point
implies that $r:=d/(p+q)$ is an integer and hence
$\alpha_1= (p+q)(r-j)$. Thus we can write
$$ C= \sum_{j=0}^r b_j(z_1^{p+q})^{r-j}(z_2^p z_3^q)^j=
 \underset{1\leq \nu \leq r}{\prod} (c z_1^{p+q} -c_\nu z_2^p z_3^q)
\quad\text{ for some } b_j,c_\nu,c \in \Bbb C.$$
Hence $X\in \Tspqc$ for $\sigma$, the identity, and
$\uc=[c_\nu]_{1\leq \nu\leq r}$. By definition $\lambda x=0\in \RF {2d-1}$
for $\forall x\in \Il d \supset w P^{d-1}+\JF d$.
Thus Th.\ref{t3-2}(3) shows
$\lambda\in \SU$. This proves the first assertion of Th.\ref{t2-1}(2).
To show the second assertion, we note that $\uc\in \Bbb P^r$
has been uniquely determined by $0\in \Dlo$. Applying the same argument to
any $t\in \Dlo$, we get a holomorphic map
$g: \Dlo \to \Bbb P^3$
defined by the condition:
$$ g(t)=\uct:=[c_{\nu,t}]_{1\leq \nu\leq r}
\text{ with } \Ft \in w_t P^{d-1} +
\underset{1\leq \nu \leq r}{\prod} (c_t z_1^{p+q} -c_{\nu,t} z_2^p z_3^q)
$$
If $\lambda\in H^2(U,\Bbb Q(2))$, then $\lt\in H^2(\Ut,\Bbb Q(2))$
for any $t\in \Delta$ and the assumption
$\lambda\not\in \JF {d-1}$ implies $\lt\not\in \JFt {d-1}$
in view of Lem.\ref{l3-1}(1). By Pr.\ref{p1-1} it implies
$\uct\in \Bbb P^r(\overline{\Bbb Q})$ and hence that $g$ is constant.
Therefore $\Dlo\subset \Tspqc \cap \Delta$ and hence
$\Dl\subset \Tspqc \cap \Delta$ by taking the closure in $\Delta$.
Finally, comparing the codimension in $\Delta$, we conclude that the last
inclusion is the equality and the proof of Th.\ref{t3-2} is complete.
\vskip 20pt

\section{Proof of key lemma}

\medbreak

In this section we prove Lem.\ref{l4-2}.
Assume that $w_t\in \sum_{i=1}^3 \Bbb C\cdot z_i$ for $\forall t\in \Dl$.
Recall that we have fixed
$\Gamma=[\gamma_1:\gamma_2:\gamma_3:\Lq] \in \dP(E_w)$
such that $0\in \Mwq$.
If $\Lq \not\in \sum_{i=1}^3 \Bbb C\cdot z_i$, (4-2) implies $\pd{F}{z_0}=0$
on $[1:0:0:0]\in \Bbb P^3$, which contradicts (3-1). Hence we have
\begin{description}
\item[(5-1)]
$\Lq\in \sum_{i=1}^3 \Bbb C\cdot z_i$
\end{description}
We may write
$$ w_t= \sum_{i=1}^3 \ait z_i\text{ and }
w=w_0= \sum_{i=1}^3 a_i z_i,$$
where $\ait$ is a holomorphic function on $\Dl$ with $a_i=\aio$.
The condition (4-3) now reads
$\sum_{i=1}^3 \gamma_i a_i z_i \in \Bbb C\cdot \sum_{i=1}^3 a_i z_i,$
which implies the condition:
\begin{description}
\item[(5-2)]
$\gamma_1 a_1 :\gamma_2 a_2 :\gamma_3 a_3  = a_1 :a_2 :a_3$.
\end{description}
The proof is now divided into some cases. Frist we suppose that we are in:

\begin{description}
\item[Case (1)] There exists $t\in \Dl$ such that
$\ait\not=0$ for $1\leq \forall i \leq 3$.
\end{description}
Without loss of generality we may suppose that $t=0$ satisfies the above
condition. (5-2) implies $\gamma_1=\gamma_2=\gamma_3$.
If $\gamma_i=0$ for $1\leq \forall i\leq 3$, $\Lq\not\in \Bbb C\cdot w$.
Then (4-2) implies $\pd{F}{z_0}\in wP^{d-1}$
so that $\pd{F}{z_0}=0$ on $[1:0:0:0]$, which contradicts (3-1).
Thus we may assume $\gamma_i=1$ for $1\leq \forall i\leq 3$.
By noting
$d F =\sum_{i=0}^3 z_i \frac{\partial F}{\partial z_i}$, (4-2) now reads:
\begin{description}
\item[(5-3)]
$ d F + (L - z_0 ) \frac{\partial F}{\partial z_0} \in w P^{d-1}.$
\end{description}

\Claim {1} \it
$\Lq\not\in \Bbb C\cdot w$.
\rm

\proof
Assume $L\in \Bbb C\cdot w$. (5-3) implies
$ d F - z_0  \frac{\partial F}{\partial z_0} \in w P^{d-1}$.
By the assumption $a_1\not=0$ we can write
$$F= w A +z_0 B + C \text{ with } B\in \Bbb C[z_0,z_2,z_3]\cap P^{d-1},\;
 C\in \Bbb C[z_2,z_3]\cap P^d.$$
Then
$\pd{F}{z_0} =w \pd {A}{z_0} +z_0 \pd{B}{z_0} +B$ and hence
$d(z_0B+C)-z_0(z_0\pd{B}{z_0}+B) =0$
by noting $wP^{d-1}\cap \Bbb C[z_0,z_2,z_3]=0$.
It implies
$C=0$ and $(d-1)B=z_0\pd{B}{z_0}$.
From the last equation we immediately deduce $B=c z_0^{d-1}$ with some
$c\in \Bbb C$. Hence $F=wA + cz_0^d$, which is singular on $\{w=A=z_0=0\}$.
It contradicts (3-1) and completes the proof of Claim 1.
\medbreak

Now choose $u\in \sum_{i=1}^3 a_i z_i$ such that $w,L,u$ are linearly
independent and write
$$F=wA + \sum_{\nu=0}^{d} L^\nu B_\nu,\text{ with }
B_\nu\in \Bbb C[z_0,u]\cap P^{d-\nu}.$$
Then the condition (5-2) implies
$$
d(\sum_{\nu=0}^{d} L^\nu B_\nu)+
(L-z_0)\sum_{\nu=0}^{d} L^\nu \pd{B_\nu}{z_0}=
\sum_{\nu=0}^{d} L^\nu (d B_\nu -z_0\pd{B_\nu}{z_0} + \pd{B_{\nu-1}}{z_0})
 \in wP^{d-1}$$
where $B_{-1}=0$ by convention. Hence we get
$ d B_\nu -z_0\pd{B_\nu}{z_0} + \pd{B_{\nu-1}}{z_0}=0$
for $0\leq \forall \nu\leq d$.
We easily solve the equations to get
$B_\nu=c (-1)^\nu \binom{d}{\nu} z_0^{d-\nu}$ for some $c\in \Bbb C$
and hence
$$F=wA + c \sum_{\nu=0}^{d} (-1)^\nu L^\nu \binom{d}{\nu} z_0^{d-\nu}
= wA + c (z_0 - L)^d.$$
The equation is singular on $\{w=A=z_0-L=0\}\subset \Bbb P^3$, which
contradicts (3-1). This completes the proof in Case (1).
\medbreak

By Case (1) we may suppose
$\Dl\subset \cup_{1\leq i\leq 3}\{t\in \Delta|\; \ait=0\}$.
Since we have shown that $\Dl$ is irreducible, we may suppose
$\act=0$ for $\forall t\in \Dl$. Now we assume that we are in:
\begin{description}
\item[Case (2)] There exists $t\in \Dl$ such that $\aat\abt\not=0$.
\end{description}

Put
$\Dlb=\{t\in \Dl|\; \aat\abt\not=0 \}$. It is a non-empty open subset of $\Dl$.
Without loss of generality we may assume $0\in \Dlb$.
Thus $a_3=0$ and $a_1a_2\not=0$.
(5-2) implies $\gamma_1=\gamma_2$. Assuming $\gamma_3\not=0$, (4-2)
implies $z_3\pd{F}{z_3}=0$ on $\{z_1=z_2=\pd{F}{z_0}=0\}$, which contradicts
(3-1).
Thus $\gamma_3=0$. If $\gamma_1=\gamma_2=0$, the same argument as in the
begining of Case (1) induces a contradiction. Thus we may assume
$\gamma_1=\gamma_2=1$. Hence (4-2) now reads:
\begin{description}
\item[(5-4)]
$ \sum_{i=1}^2 z_i\frac{\partial F}{\partial z_i}+
 L \frac{\partial F}{\partial z_0} \in wP^{d-1}.$
\end{description}

\Claim {2}
$\Lq\in \sum_{i=1}^2 \Bbb C\cdot z_i$ and $L\not\in \Bbb C\cdot w$.
\rm

\proof
Assume $\Lq \not\in \sum_{i=1}^2 \Bbb C\cdot z_i$. By (5-1) we may suppose
$\Lq=z_3+l_1z_1+l_2 z_2$. Then (5-4) implies
$\pd{F}{z_0}\in <z_1,z_2>$, which contradicts (3-1).
The proof of the second assertion is smilar to that of Claim 1 and omitted.
\medbreak

Noting $\Bbb C[z_0,z_1,z_2,z_3]=\Bbb C[z_0,w,L,z_3]$, we may write
$$ F= wA + \sum_{\mu=0}^d z_3^\mu G_\mu
\text{ with } A\in P^{d-1},\; G_\mu\in \Bbb C[z_0,L]\cap P^{d-\mu}.$$
Noting $\sum_{i=1}^2 z_i \pd{w}{z_i}=w$, (5-4) implies
$$ \sum_{\mu=0}^d z_3^\mu
(\sum_{i=1}^2 z_i\frac{\partial G_\mu}{\partial z_i}+
L \frac{\partial G_\mu}{\partial z_0})  \in wP^{d-1}.$$
Noting $(d-\mu)G_\mu=\sum_{i=0}^3\pd{G_\mu}{z_i}$ and $\pd{G_\mu}{z_3}=0$,
we get
$$ 0= \sum_{i=1}^2 z_i\frac{\partial G_\mu}{\partial z_i}+
L \frac{\partial G_\mu}{\partial z_0}
=(d-\mu)G_\mu+(L-z_0)\pd{G_\mu}{z_0} \text{ for }\forall 1\leq \mu\leq d.$$
We solve the last equation in the same manner as Case (1) to get
$G_\mu= b_\mu (L-z_0)^{d-\mu}$ with $b_\mu\in \Bbb C$ and hence
\begin{itemize}
\item[\bf(5-5)]
$F= wA + \sum_{\mu=0}^d b_\mu z_3^\mu (L-z_0)^{d-\mu}.$
\end{itemize}

\Claim {3} \it Put
$\eF= A +\sum_{i=1}^2 z_i\pd{A}{z_i} + L \pd{A}{z_0}$.
\begin{itemize}
\item[$(1)$]
$\phi(\eF)=\frac{z_3}{w} d(\frac{z_0- L}{z_3})\wedge dlog{\frac{z_1}{z_2}}$.
\item[$(2)$]
$\Bbb C\cdot \wF\oplus\Bbb C\cdot \eF =
\{y\in P^{d-1}|\; yx=0\in \RF {2d-1}\text{ for }\forall
x\in wP^{d-1}+\JF d\}$  (cf. Lem.\ref{l3-1}).
\end{itemize}
\rm
\medbreak
Claim 3(1) is easily proven by using (3-\ref{i3-3}) and Claim 3(2)
is proven by the same argument as the proof of Th.\ref{t3-2}(3).
We omit the details.
\medbreak

By Claim 3 $\lt\in \HUt 20$, the flat translation of
$\lambda$ for $t\in \Dl$, is written as
$$ \lt = \fat \et + \fbt \wt \text{ for } t\in \Dlb.$$
Here $\fat$ and $\fbt$ are holomorphic functions on $\Dlb$ and
$$\wt= \dlog{\frac{z_2}{z_1}}\wedge \dlog{\frac{z_3}{z_1}},\quad
\et=\frac{z_3}{w_t} d(\frac{z_0- L_t}{z_3})\wedge dlog{\frac{z_1}{z_2}},$$
where $w_t$ is as in the begining of this section and
$$
\Ft=w_t A_t + \sum_{\mu=0}^d b_{\mu,t} z_3^\mu (L_t-z_0)^{d-\mu},
\quad L_t= l_1(t) z_1 + l_2(t) z_2
$$
is the defining equation of $\Xt$ such as (5-5), which varies holomorphically
with $t\in \Dlb$.
Recalling $Y=\cup_{1\leq j\leq 3} Y_j$ with $Y_j=\{z_j=0\}\subset \Bbb P^3$,
write
$$
Z_t=\Xt\cap Y \supset Z_{3t}=\Xt\cap Y_3 \supset
V_t=Z_{3t} \cap ( Y_1 \cup Y_2) \supset
S_t=Z_{3t} \cap Y_2.$$
We consider the composite of the residue maps
$$\theta_t \scs \HUt 20=H^0(\Xt,\WXZt 2) \rmapo{Res_{Z_{3t}}}
H^0(Z_{3t},\Omega_{Z_{3t}}^1(\log V_t)) \rmapo{Res_{S_t}}
\Bbb C^{S_t} \isom \Bbb C^d,$$
where the last isomorphism is obtained by choosing
$\epsilon_t: \{1,2,\dots,d\} \isom S_t$, an isomorphism of local systems of
sets over $\Delta$. Since $\lt$ is flat, we get the condition:
\begin{description}
\item[(5-6)]
$\theta_t(\lt)\in \Bbb C^d$ is constant with $t\in \Dlb$.
\end{description}
We shall show that the condition induces a contradiction, which completes
the proof of Lem.\ref{l4-2} in Case (2).
In order to calculate $\theta_t(\lt)$ we introduce some notations.
Let $\bA=\{z_2=z_3=0\}-\{[0:1:0:0]\}\subset \Bbb P^3$ be identified with
$\Bbb C$ via $[z_0:z_1]\to z_0/z_1$. Let
$$ \Sigma=\{(s_1,\dots,s_d)|\; s_\nu\in \bA,\; s_\nu\not=s_\mu
\text{ for } 1\leq \nu\not=\mu\leq d \}.$$
We define a holomorphic map
$$ \pi \scs \Delta \to \Sigma;\; t \to (\snt)_{1\leq \nu\leq d}
\text{ with } \epsilon_t(\nu)=\snt.$$
Now an easy residue calculation shows
$$ \theta_t(\wt)=(1,\dots,1),\quad
   \theta_t(\et)=(\frac{l_1(t)- \snt}{a_1(t)})_{1\leq \nu\leq d}$$
and hence
$$ \theta_t(\lt)=(p(t)\snt + q(t))_{1\leq \nu\leq d}\text{ with }
p(t)=-\frac{\fat}{a_1(t)},\; q(t)=\fat\frac{l_1(t)}{\aat}+ \fbt.$$
Therefore (5-6) implies that for $\forall \partial \in T_0(\Dl)$, we have
$$ 0= \partial(p(t)\snt + q(t))= p(0) \partial\snt + \sno \partial p(t) +
\partial q(t)\text{ for } 1\leq \forall \nu\leq d.$$
Letting $\pi_*: T_0(\Delta) \to T_{\pi(0)}(\Sigma)\cong \Bbb C^d$
be the differential of $\pi$, we get
$$ p(0)\cdot \pi_*(\partial) =
p(0)\cdot(\partial\snt)_{1\leq \nu\leq d}=
-\partial p(t)\cdot (\sno)_{1\leq \nu\leq d}+\partial q(t)\cdot (1,\dots,1).$$
Since $\lambda\not\in \Bbb C\cdot \wF$,
$\fao\not=0$ and hence $p(0)\not=0$.
Thus it implies $\dim(\pi_*(T_0(\Dl)))\leq 2$.
Therefore we get a contradiction if we show the following.

\Claim {4} \it
$\dim(\pi_*(T_0(\Dl)))\geq d$.
\rm

\proof
Let $Q=\Bbb C[z_0,z_1]$ and $Q^l= P^l \cap Q$ for an integer $l\geq 0$.
For $G\in P$ write $\overline{G}=G \mod <z_2,z_3> \in Q$. Consider the morphism
$$ \rho \scs \Sigma \to N:=\dP(Q^d);\; \us=(s_\nu)_{1\leq \nu\leq d} \to
[F_{\us}]\text{ with } F_{\us}=\underset{1\leq \nu\leq d}{\prod}
(z_0 -s_\nu z_1).$$
It is finite etale and induces an isomorphism on the tangent spaces.
Hence it suffices to show Claim 4 by replacing $\pi$ with
$\tilde{\pi}:=\rho\circ\pi$. We have
$ \tilde{\pi}(t)=[\overline{F}_t]$
and we have the commutative diagram
$$
  \begin{matrix}
   P^d &\overset{\mod <z_2,z_3>}{\longrightarrow}& Q^d \\
   \downarrow\rlap{$\psi$} && \downarrow\rlap{$\psi'$} \\
   T_0(\Delta) & \overset{\tilde{\pi}}{\longrightarrow} &
   T_{\tilde{\pi}(0)}(N) \\
   \end{matrix}
  $$
where $\psi'$ is defined in the same say as $\psi$ in (3-1) and
$\Ker(\psi')=\Bbb C\cdot \overline{F}$. We have shown that
$\psi^{-1}(T_0(\Dl))=\Il d \supset w P^{d-1} +\JF d$. Hence
$\tilde{\pi}_*(T_0(\Dl)\supset \psi'(z_1 Q^{d-1} +\Bbb C\cdot\overline{F})$.
Noting $F\not\in <z_1,z_2,z_3>$ so that $\overline{F}\not\in z_1 Q^{d-1}$,
this implies
$$ \dim(\tilde{\pi}_*(T_0(\Dl))\geq \dim z_1 Q^{d-1}=d.$$
This completes the proof of Claim 4.
\medbreak

By Case (2) we may assume now that we are in:
\begin{description}
\item[Case (3)]
$\abt=\act=0$ for $\forall t\in \Dl$.
\end{description}
In this case we may assume $w=z_1$. We have
$$I_\lambda \supset I:=<z_1> + J_F =
<z_1,\pd{F}{z_0},z_2\pd{F}{z_2},z_3\pd{F}{z_3}>$$
so that $I$ is complete intersection of degree $(1,d-1,d,d)$.
Hence $I=I_\lambda$ and $\Il d=z_1 P^{d-1} +\JF d$.
As before we can show the following.

\Claim {5} \it Put $\kF=\pd{F}{z_1}$.
\begin{itemize}
\item[$(1)$]
$\phi(\kF)=\frac{z_0}{z_1} \dlog{\frac{z_2}{z_0}}
\wedge dlog{\frac{z_3}{z_0}}$.
\item[$(2)$]
$\Bbb C\cdot \wF\oplus\Bbb C\cdot \kF =
\{y\in P^{d-1}|\; yx=0\in \RF {2d-1}\text{ for }\forall
x\in z_1P^{d-1}+\JF d\}$.
\end{itemize}
\rm
\medbreak

As before Claim 5 implies
$$ \lt = \fat \kt + \fbt \wt \text{ for } t\in \Dl,$$
where $\fat$, $\fbt$ and $\wt$ are as before and
$$ \kt=\frac{z_0}{z_1} \dlog{\frac{z_2}{z_0}}\wedge dlog{\frac{z_3}{z_0}}
\in \HUt 20. $$
An easy residue caluculation shows
$\theta_t(\lt)=(\fat \snt + \fbt)_{1\leq \nu\leq d}$
and the same argument as Case (2) induces a contradiction.
This completes the proof of Lem.\ref{l4-2}.
\vskip 20pt

\section{Injectivity of regulator map}
\medbreak

In this section we prove Th.\ref{0-6}.
Fix $0\in M$ and let the notaion be as in the begining of \S1.
By Lem.\ref{l1-1}, if $0\in M\setminus \MNL$, we have
$$ F^0H^2(U,\Bbb Q(2)) =\Bbb Q\cdot \regU 2 (\aU)
\text{ with }
\aU =\{\frac{z_2}{z_1},\frac{z_3}{z_1}\}\in CH^2(U,2).$$
By [AS2], Th.(6-1), it implies that the kernel of the composite map
$$ CH^1(Z,1)\otimes\Bbb Q \to CH^2(X,1)\otimes\Bbb Q
\rmapo{\oX} H^3_D(X,\Bbb Q(2))$$
is generated by
$$ \partial_{U}(\aU)=\delta:=
((\frac{z_3}{z_2})_{|Z_1},(\frac{z_1}{z_3})_{|Z_2},(\frac{z_2}{z_1})_{|Z_2})
\in CH^1(Z,1),$$
where
$\partial_{U}:CH^2(U,2) \to CH^1(Z,1)$.

\Claim {1} \it Write $\Lambda= \bigoplus_{1\leq j\leq 3} \Bbb C(Z_j)^*$.
\begin{itemize}
\item[$(1)$]
Assume $0\in \Ta$ and that $X$ is defined by an equation as Def.\ref{d0-3}:
$$ F=w A + z_1 z_2 B + c_1 z_1^d + c_2 z_2^d.$$
Then the following elements are linearly independent in $\Lambda\otimes\Bbb Q$;
$$
\delta,
\; ((\frac{z_2}{w})_{|Z_1},(\frac{w}{z_1})_{|Z_2},1).$$
\item[$(2)$]
Assume $0\in \Ta\cap\Tb$ and that $X$ is defined
by an equation as Def.\ref{d0-3}:
\begin{align*}
 F & =w A + z_1 z_2 B + c_1 z_1^d + c_2 z_2^d \\
   & =v A + z_2 z_3 B' + c_2 z_2^d + c_3 z_3^d
\end{align*}
Then the following elements are linearly independent in $\Lambda\otimes\Bbb Q$;
$$
\delta,
\; ((\frac{z_2}{w})_{|Z_1},(\frac{w}{z_1})_{|Z_2},1),
\; (1,(\frac{z_3}{v})_{|Z_2},(\frac{v}{z_2})_{|Z_3} ).
$$
\item[$(3)$]
Assume $0\in \Ta\cap\Tb$ and that $X$ is defined
by an equation as Def.\ref{d0-3}:
\begin{align*}
 F & =w A + z_1 z_2 B + c_1 z_1^d + c_2 z_2^d \\
   & =v A + z_2 z_3 B' + c_2 z_2^d + c_3 z_3^d \\
   & =u A + z_3 z_1 B'' + c_3 z_3^d + c_1 z_1^d
\end{align*}
Then the following elements are linearly independent in $\Lambda\otimes\Bbb Q$;
$$
\delta,
\; ((\frac{z_2}{w})_{|Z_1},(\frac{w}{z_1})_{|Z_2},1),
\; (1,(\frac{z_3}{v})_{|Z_2},(\frac{v}{z_2})_{|Z_3} ),
\; ((\frac{u}{z_3})_{|Z_1},1,(\frac{z_1}{u})_{|Z_3}).
$$
\end{itemize}
\rm

\proof
We only show Claim 1(3). The other are easier and left to the readers.
Assume the contrary. Then there are integers $e,l,m,n$ not all zero such that
\begin{align*}
&(\frac{z_2}{w})^l (\frac{u}{z_3})^n (\frac{z_3}{z_2})^e \equiv 1\mod z_1,\\
&(\frac{w}{z_1})^l (\frac{z_3}{v})^m (\frac{z_1}{z_3})^e \equiv 1\mod z_2,\\
&(\frac{v}{z_2})^m (\frac{z_1}{u})^n (\frac{z_2}{z_1})^e \equiv 1\mod z_3.
\end{align*}
We note $u,v,w\not\in \sum_{1\leq j\leq 3}\Bbb C\cdot z_j$ since
otherwise it would contradicts (3-1).
Hence the condition implies $l=m=n=e$ and
$u,v,w$ coincides up to non-zero constant. Thus we get
$$ F  \equiv
   z_1 z_2 B + c_1 z_1^d + c_2 z_2^d \equiv
   z_2 z_3 B' + c_2 z_2^d + c_3 z_3^d \equiv
   z_3 z_1 B'' + c_3 z_3^d + c_1 z_1^d\mod w,$$
which is absurd. This completes the proof of Claim 1.
\medbreak

By Claim 1, the proof of Th.\ref{0-6} is complete if we show that
$\Ta\not\subset\MNL$ (resp. $\Ta\cap\Tb\not\subset\MNL$,
resp. $\Ta\cap\Tb\cap\Tc\not\subset\MNL$) if $d\geq 4$ (resp. $d\geq 6$,
resp. $d\geq 10$). Indeed we have
$$ \codim_M(\Ta)=\binom{d+3}{3}-
\left(\binom{d+2}{3}+\binom{d}{2}+2\right)=2d-1.$$
One note that $\Ta\cap \Tb\cap \Tc\not=\emptyset$ since
the Fermat surface $z_0^d+z_1^d+z_2^d+z_3^d=0$ belongs to it.
Hence, for any irreducible component $T$ of $\Ta\cap\Tb$
(resp. $\Ta\cap\Tb\cap\Tc$),
$\codim_M(T)\leq 2(2d-1)$ (resp. $\codim_M(T)\leq 3(2d-1)$).
By Th.\ref{0-4}(1) it suffices to check $\binom{d+2}{2}-5$ is greater than
$2d-1$ (resp. $2(2d-1)$, resp. $3(2d-1)$) if $d\geq 4$ (resp. $d\geq 6$,
resp. $d\geq 10$). This completes the proof of Th.\ref{0-6}.

\end{document}